\documentclass[reqno,12pt]{article}
\usepackage{amsmath}
\usepackage{amsthm}
\usepackage{amssymb}
\usepackage{enumerate}
\usepackage{authblk}
\usepackage[T1]{fontenc}
\usepackage{float}
\usepackage{tikz-network}
\usepackage{pdfpages}

  \providecommand{\definitionname}{Definition}
\newtheorem{theorem}{Theorem}

\newtheorem{proposition}[theorem]{Proposition}

\begin{document}

\title{Online size Ramsey number for $C_4$ and $P_6$}
\author{Mateusz Litka}
\affil{Adam Mickiewicz University in Poznań}

\maketitle

\begin{abstract}

In this paper we consider a game played on the edge set of the infinite clique $K_\mathbb{N}$ by two players, \textit{Builder} and \textit{Painter}. In each round of the game, Builder \textit{chooses} an edge and Painter \textit{colors} it red or blue. Builder wins when Painter creates a red copy of $G$ or a blue copy of $H$, for some fixed graphs $G$ and $H$. Builder wants to win in as few rounds as possible, and Painter wants to delay Builder for as many rounds as possible.\newline \indent
\textit{The online size Ramsey number} $\tilde{r}(G,H)$, is the minimum number of rounds within which Builder can win, assuming both players play optimally.\newline \indent
So far it has been proven by Dybizbański, Dzido and Zakrzewska that $11\leq\tilde{r}(C_4,P_6)\leq13$ \cite{Dzido}. In this paper, we refine this result and show the exact value, namely we will present the Theorem that $\tilde{r}(C_4,P_6)=11$, with the details of the proof.
\end{abstract}

\section{Introduction}

In this paper, we consider the following generalization of Ramsey numbers, which was introduced by Beck \cite{Beck}. Consider a game played on the edge set of the infinite clique $K_\mathbb{N}$ by two players, \textit{Builder} and \textit{Painter}. In each round of the game, Builder \textit{chooses} an edge and Painter \textit{colors} it red or blue. Builder wins when Painter creates a red copy of $G$ or a blue copy of $H$, for some fixed graphs $G$ and $H$. Builder wants to win in as few rounds as possible, and Painter wants to delay Builder for as many rounds as possible.\newline \indent
By red-blue \textit{edge coloring} of a graph $G$ we mean a function defined on the set of edges $E(G)$ that assigns one of the colors - red or blue - to each of the edges. We say that a graph $G$ is \textit{colored} if each of its edges is either blue or red. A graph is \textit{red} if all its edges are red. A graph is \textit{blue} if all its edges are blue. Such graphs are called \textit{monochromatic} (1-colored).\newline \indent
We say that Builder \textit{forces} a red edge if, after Builder selects that edge, then Painter coloring it blue will create a blue copy of $H$. Similarly, we say Builder \textit{forces} a blue edge if, after Builder selects that edge, then Painter coloring it red will create a red copy of $G$.\newline \indent
We say that Builder \textit{forces} a red copy of $G$ or a blue copy of $H$ in a given number of rounds $t$ if it has a strategy that ensures that after at most $t$ rounds the edge selected by Builder forces both red and blue edge at the same time. Then Painter coloring it blue will create a blue $H$, and coloring it red will create a red $G$. This prevents Painter from making another non-losing move.\newline \indent
\textit{The online size Ramsey number} of two graphs $G$ and $H$, denoted by $\tilde{r}(G,H)$, is the minimum number of rounds in which Builder forces a red copy of $G$ or a blue copy of $H$, assuming both Builder and Painter are playing optimally.\newline \indent
We will call this game \textit{$\tilde{R}(G,H)$-game}. First note that $\tilde{r}(G,H)=\tilde{r}(H,G)$. One can also consider a symmetric version of the numbers defined above. If $H=G$, then instead of $\tilde{r}(G,G)$ we write $\tilde{r}(G)$.

\section{Main result}

As mentioned above, online size Ramsey number was introduced by Beck \cite{Beck}. The exact values of $\tilde{r}(G,H)$ are known for several graphs $G$ and $H$. In \cite{Rucinski}, Kurek and Ruciński considered the case where $G$ and $H$ are cliques, but except for the trivial $\tilde{r}(K_2,K_k)={k \choose 2}$, they were able to determine only one more value, namely $\tilde{r}(K_3,K_3)=8$. So far, only one more value is known, namely $\tilde{r}(K_3,K_4)=17$, which was obtained by Prałat \cite{Pralat4} with computer support. Grytczuk, Kierstead and Prałat \cite{Grytczuk} gave exact values for several short paths and showed that $\tilde{r}(P_6)=10$. Prałat also showed that $\tilde{r}(P_7)=12$, $\tilde{r}(P_8)=15$ \cite{Pralat2} and $\tilde{r}(P_9)=17$ \cite{Pralat3}.\newline \indent
We now turn to the subject of this paper, i.e. to the online size Ramsey numbers of the form $\tilde{r}(C_4,P_k)$. Values $\tilde{r}(C_4,P_3)=6$ and $\tilde{r}(C_4,P_4)=8$ were obtained by Cyman, Dzido, Lapinskas and Lo \cite{Lo}. Then in \cite{Dzido} Dybizbański, Dzido and Zakrzewska proved that $\tilde{r}(C_4,P_5)=9$. Also in \cite{Dzido} were obtained bounds for longer paths.
\begin{proposition}(\cite{Dzido})\label{bounds}
	For $k\geq6$ we have $$3k-5\geq\tilde{r}(C_4,P_k)\geq \left\{ \begin{array}{ll}
	2k-1 & \textrm{if $k=6,7;$}\\
	2k-2 & \textrm{if $k\geq8. $}
	\end{array} \right. $$
\end{proposition}
Thus, from the results obtained so far, it follows that $13\geq\tilde{r}(C_4,P_6)\geq11$. In this paper, we prove the following Theorem.
\begin{theorem}\label{main}
	$\tilde{r}(C_4,P_6)=11$.
\end{theorem}

\section{Proof of Theorem \ref{main}}

\begin{figure}[H]
	\centering
	\includegraphics[width=0.47\linewidth]{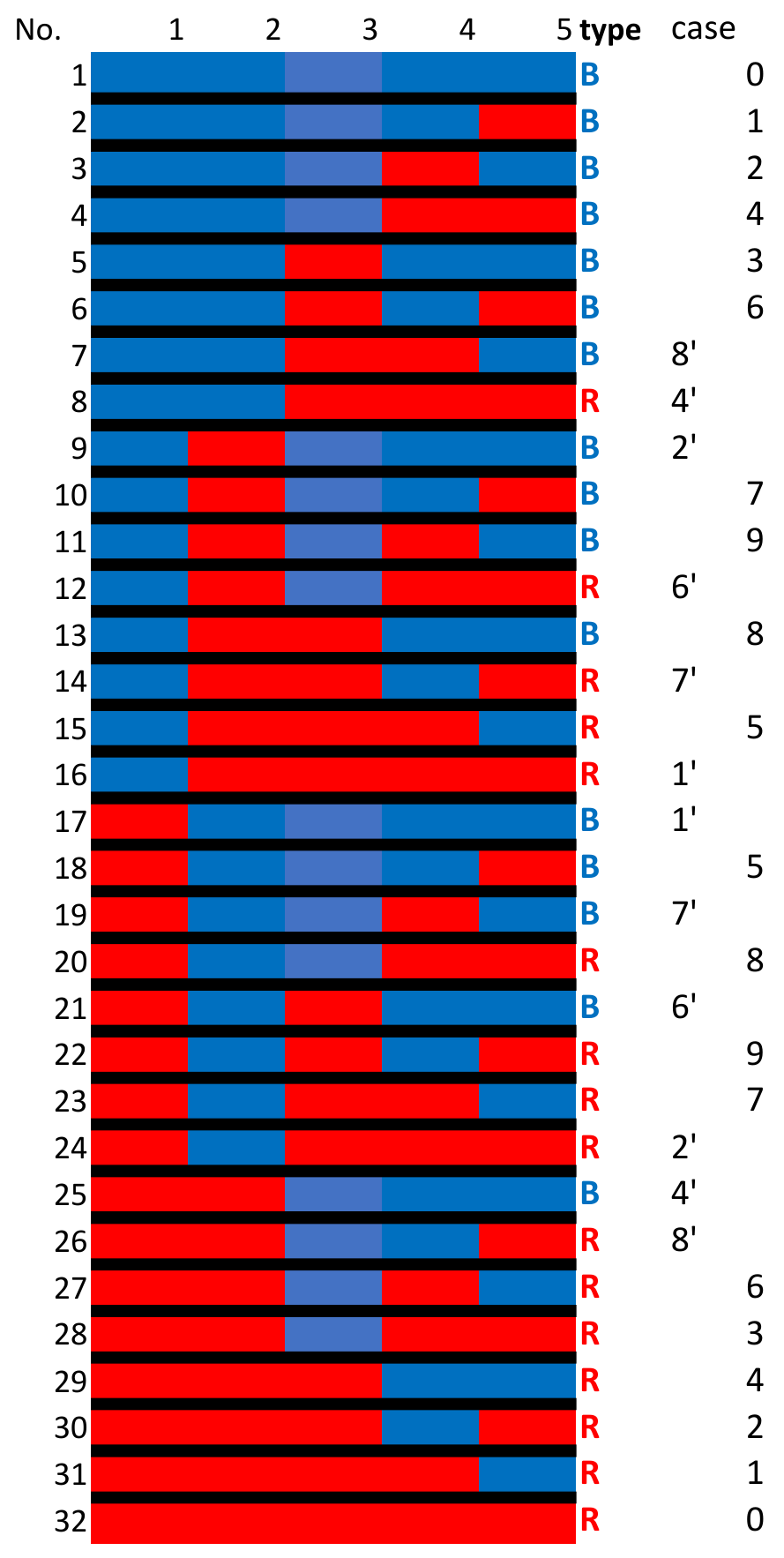}
	\caption[1]{First five rounds.}
	\label{1}
\end{figure}

By Proposition \ref{bounds} $11\leq\tilde{r}(C_4,P_6)\leq13$. Therefore, it suffices to prove that Builder can win the $\tilde{R}(C_4,P_6)$-game in $11$ rounds. We will present a strategy for Builder that forces Painter to create a red copy of $C_4$ or a blue copy of $P_6$ in at most $11$ rounds. The upper bound $(\tilde{r}(C_4,P_6)\leq11)$ is obtained by analyzing individual cases, presented in the table in Figure \ref{1} (first five rounds) and in the graphs below (last six rounds).
\subsection*{First five rounds}
In the first five rounds, Builder constructs the path $P_6$. After selecting five edges of the path $P_6$ in sequence, $2^5=32$ of color patterns are possible. They are presented and numbered in the table in Figure \ref{1}. Since we only have two colors (blue and red) and $P_6$ has an odd number of edges ($5$) to color, thus we will always get one of two types: $B$ (from blue) or $R$ (from red), where we have more edges colored blue or red, respectively.\newline \indent
Basically, one of the $19$ possible color patterns appears (up to the symmetry), as indicated in the table by the appropriate case number (numbers from $0$ to $9$). Cases that are symmetrical to those already considered are marked by adding the symbol $'$ (apostrophe). Note that we have discarded the $B0$ case, because in this case Builder would win after the first $5$ rounds.
\subsection*{Last six rounds}
The rest of the result is obtained by analyzing the individual subcases of the last six rounds, which are shown in the graphs below.\newline\newline
\noindent B1(a)
\begin{tikzpicture}[thick,scale=0.8, every node/.style={scale=0.8}]
\Vertex[x=0,y=0,IdAsLabel]{A}
\Vertex[x=0,y=2,IdAsLabel]{B}
\Vertex[x=1,y=4,IdAsLabel]{C}
\Vertex[x=3,y=4,IdAsLabel]{D}
\Vertex[x=4,y=2,IdAsLabel]{E}
\Vertex[x=4,y=0,IdAsLabel]{F}
\Vertex[x=2,y=2,IdAsLabel]{G}
\Edge[color=blue,label=1,position=left](A)(B)
\Edge[color=blue,label=2,position=above left](B)(C)
\Edge[color=blue,label=3,position=above](C)(D)
\Edge[color=blue,label=4,position=above right](D)(E)
\Edge[color=red,label=5,position=right](E)(F)
\Edge[color=red,label=6,position=below](A)(F)
\Edge[color=red,label=7,position=above left](A)(G)
\Edge[color=blue,label=8,position=above,bend=5](G)(E)
\Edge[color=red,label=8,position=below,bend=-5](G)(E)
\end{tikzpicture}
\\
B2(a)
\begin{tikzpicture}[thick,scale=0.8, every node/.style={scale=0.8}]
\Vertex[x=0,y=0,IdAsLabel]{A}
\Vertex[x=0,y=2,IdAsLabel]{B}
\Vertex[x=1,y=4,IdAsLabel]{C}
\Vertex[x=3,y=4,IdAsLabel]{D}
\Vertex[x=4,y=2,IdAsLabel]{E}
\Vertex[x=4,y=0,IdAsLabel]{F}
\Edge[color=blue,label=1,position=left](A)(B)
\Edge[color=blue,label=2,position=above left](B)(C)
\Edge[color=blue,label=3,position=above](C)(D)
\Edge[color=red,label=4,position=above right](D)(E)
\Edge[color=blue,label=5,position=right](E)(F)
\Edge[color=red,label=6,position=above left](A)(E)
\Edge[color=red,label=7,position=left,bend=-15](D)(F)
\Edge[color=blue,label=8,position=above,bend=5](A)(F)
\Edge[color=red,label=8,position=below,bend=-5](A)(F)
\end{tikzpicture}
\\
B3(a)
\begin{tikzpicture}[thick,scale=0.8, every node/.style={scale=0.8}]
\Vertex[x=0,y=0,IdAsLabel]{A}
\Vertex[x=0,y=2,IdAsLabel]{B}
\Vertex[x=1,y=4,IdAsLabel]{C}
\Vertex[x=3,y=4,IdAsLabel]{D}
\Vertex[x=4,y=2,IdAsLabel]{E}
\Vertex[x=4,y=0,IdAsLabel]{F}
\Edge[color=blue,label=1,position=left](A)(B)
\Edge[color=blue,label=2,position=above left](B)(C)
\Edge[color=red,label=3,position=above](C)(D)
\Edge[color=blue,label=4,position=above right](D)(E)
\Edge[color=blue,label=5,position=right](E)(F)
\Edge[color=red,label=6,position=above left](A)(D)
\Edge[color=red,label=7,position=above right](C)(F)
\Edge[color=blue,label=8,position=above,bend=5](A)(F)
\Edge[color=red,label=8,position=below,bend=-5](A)(F)
\end{tikzpicture}
\\
B4(a)
\begin{tikzpicture}[thick,scale=0.8, every node/.style={scale=0.8}]
\Vertex[x=0,y=0,IdAsLabel]{A}
\Vertex[x=0,y=2,IdAsLabel]{B}
\Vertex[x=1,y=4,IdAsLabel]{C}
\Vertex[x=3,y=4,IdAsLabel]{D}
\Vertex[x=4,y=2,IdAsLabel]{E}
\Vertex[x=4,y=0,IdAsLabel]{F}
\Vertex[x=2,y=2,IdAsLabel]{G}
\Edge[color=blue,label=1,position=left](A)(B)
\Edge[color=blue,label=2,position=above left](B)(C)
\Edge[color=blue,label=3,position=above](C)(D)
\Edge[color=red,label=4,position=above right](D)(E)
\Edge[color=red,label=5,position=right](E)(F)
\Edge[color=blue,label=$6_1$,position=below](A)(F)
\Edge[color=red,label=7,position=above left](D)(G)
\Edge[color=blue,label=8,position=above right,bend=-5](F)(G)
\Edge[color=red,label=8,position=below left,bend=5](F)(G)
\end{tikzpicture}
B4(b)
\begin{tikzpicture}[thick,scale=0.8, every node/.style={scale=0.8}]
\Vertex[x=0,y=0,IdAsLabel]{A}
\Vertex[x=0,y=2,IdAsLabel]{B}
\Vertex[x=1,y=4,IdAsLabel]{C}
\Vertex[x=3,y=4,IdAsLabel]{D}
\Vertex[x=4,y=2,IdAsLabel]{E}
\Vertex[x=4,y=0,IdAsLabel]{F}
\Vertex[x=2,y=2,IdAsLabel]{G}
\Edge[color=blue,label=1,position=left](A)(B)
\Edge[color=blue,label=2,position=above left](B)(C)
\Edge[color=blue,label=3,position=above](C)(D)
\Edge[color=red,label=4,position=above right](D)(E)
\Edge[color=red,label=5,position=right](E)(F)
\Edge[color=red,label=$6_2$,position=below](A)(F)
\Edge[color=blue,label=$7_1$,position=below left,bend=-15](D)(F)
\Edge[color=red,label=8,position=above left](A)(E)
\Edge[color=red,label=9,position=above right](F)(G)
\Edge[color=blue,label=10,position=below right,bend=-5](A)(G)
\Edge[color=red,label=10,position=above left,bend=5](A)(G)
\end{tikzpicture}
B4(c)
\begin{tikzpicture}[thick,scale=0.8, every node/.style={scale=0.8}]
\Vertex[x=0,y=0,IdAsLabel]{A}
\Vertex[x=0,y=2,IdAsLabel]{B}
\Vertex[x=1,y=4,IdAsLabel]{C}
\Vertex[x=3,y=4,IdAsLabel]{D}
\Vertex[x=4,y=2,IdAsLabel]{E}
\Vertex[x=4,y=0,IdAsLabel]{F}
\Vertex[x=2,y=2,IdAsLabel]{G}
\Edge[color=blue,label=1,position=left](A)(B)
\Edge[color=blue,label=2,position=above left](B)(C)
\Edge[color=blue,label=3,position=above](C)(D)
\Edge[color=red,label=4,position=above right](D)(E)
\Edge[color=red,label=5,position=right](E)(F)
\Edge[color=red,label=$6_2$,position=below](A)(F)
\Edge[color=red,label=$7_2$,position=below left,bend=-15](D)(F)
\Edge[color=blue,label=8,position=above left](A)(E)
\Edge[color=red,label=9,position=above](E)(G)
\Edge[color=blue,label=10,position=above left,bend=-5](D)(G)
\Edge[color=red,label=10,position=below right,bend=5](D)(G)
\end{tikzpicture}
\\
B5(a)
\begin{tikzpicture}[thick,scale=0.8, every node/.style={scale=0.8}]
\Vertex[x=0,y=0,IdAsLabel]{A}
\Vertex[x=0,y=2,IdAsLabel]{B}
\Vertex[x=1,y=4,IdAsLabel]{C}
\Vertex[x=3,y=4,IdAsLabel]{D}
\Vertex[x=4,y=2,IdAsLabel]{E}
\Vertex[x=4,y=0,IdAsLabel]{F}
\Vertex[x=2,y=2,IdAsLabel]{G}
\Edge[color=red,label=1,position=left](A)(B)
\Edge[color=blue,label=2,position=above left](B)(C)
\Edge[color=blue,label=3,position=above](C)(D)
\Edge[color=blue,label=4,position=above right](D)(E)
\Edge[color=red,label=5,position=right](E)(F)
\Edge[color=blue,label=$6_1$,position=above](A)(E)
\Edge[color=red,label=7,position=below](B)(F)
\Edge[color=red,label=8,position=below](A)(F)
\Edge[color=red,label=9,position=above left](A)(G)
\Edge[color=blue,label=10,position=above,bend=5](B)(G)
\Edge[color=red,label=10,position=below,bend=-5](B)(G)
\end{tikzpicture}
B5(b)
\begin{tikzpicture}[thick,scale=0.8, every node/.style={scale=0.8}]
\Vertex[x=0,y=0,IdAsLabel]{A}
\Vertex[x=0,y=2,IdAsLabel]{B}
\Vertex[x=1,y=4,IdAsLabel]{C}
\Vertex[x=3,y=4,IdAsLabel]{D}
\Vertex[x=4,y=2,IdAsLabel]{E}
\Vertex[x=4,y=0,IdAsLabel]{F}
\Vertex[x=2,y=2,IdAsLabel]{G}
\Edge[color=red,label=1,position=left](A)(B)
\Edge[color=blue,label=2,position=above left](B)(C)
\Edge[color=blue,label=3,position=above](C)(D)
\Edge[color=blue,label=4,position=above right](D)(E)
\Edge[color=red,label=5,position=right](E)(F)
\Edge[color=red,label=$6_2$,position=above](A)(E)
\Edge[color=blue,label=7,position=below](B)(F)
\Edge[color=red,label=8,position=below](A)(F)
\Edge[color=red,label=9,position=above right](F)(G)
\Edge[color=blue,label=10,position=above,bend=-5](E)(G)
\Edge[color=red,label=10,position=below,bend=5](E)(G)
\end{tikzpicture}
\\
B6(a)
\begin{tikzpicture}[thick,scale=0.8, every node/.style={scale=0.8}]
\Vertex[x=0,y=0,IdAsLabel]{A}
\Vertex[x=0,y=2,IdAsLabel]{B}
\Vertex[x=1,y=4,IdAsLabel]{C}
\Vertex[x=3,y=4,IdAsLabel]{D}
\Vertex[x=4,y=2,IdAsLabel]{E}
\Vertex[x=4,y=0,IdAsLabel]{F}
\Vertex[x=1,y=0,IdAsLabel]{G}
\Edge[color=blue,label=1,position=left](A)(B)
\Edge[color=blue,label=2,position=above left](B)(C)
\Edge[color=red,label=3,position=above](C)(D)
\Edge[color=blue,label=4,position=above right](D)(E)
\Edge[color=red,label=5,position=right](E)(F)
\Edge[color=blue,label=$6_1$,position=below right](A)(D)
\Edge[color=red,label=7,position=left](C)(G)
\Edge[color=red,label=8,position=above right](C)(F)
\Edge[color=blue,label=9,position=above left,bend=5](G)(E)
\Edge[color=red,label=9,position=below right,bend=-5](G)(E)
\end{tikzpicture}
B6(b)
\begin{tikzpicture}[thick,scale=0.8, every node/.style={scale=0.8}]
\Vertex[x=0,y=0,IdAsLabel]{A}
\Vertex[x=0,y=2,IdAsLabel]{B}
\Vertex[x=1,y=4,IdAsLabel]{C}
\Vertex[x=3,y=4,IdAsLabel]{D}
\Vertex[x=4,y=2,IdAsLabel]{E}
\Vertex[x=4,y=0,IdAsLabel]{F}
\Edge[color=blue,label=1,position=left](A)(B)
\Edge[color=blue,label=2,position=above left](B)(C)
\Edge[color=red,label=3,position=above](C)(D)
\Edge[color=blue,label=4,position=above right](D)(E)
\Edge[color=red,label=5,position=right](E)(F)
\Edge[color=red,label=$6_2$,position=below right](A)(D)
\Edge[color=blue,label=$7_1$,position=below](A)(F)
\Edge[color=red,label=8,position=left,bend=15](F)(D)
\Edge[color=blue,label=9,position=above right,bend=5](C)(E)
\Edge[color=red,label=9,position=below left,bend=-5](C)(E)
\end{tikzpicture}
B6(c)
\begin{tikzpicture}[thick,scale=0.8, every node/.style={scale=0.8}]
\Vertex[x=0,y=0,IdAsLabel]{A}
\Vertex[x=0,y=2,IdAsLabel]{B}
\Vertex[x=1,y=4,IdAsLabel]{C}
\Vertex[x=3,y=4,IdAsLabel]{D}
\Vertex[x=4,y=2,IdAsLabel]{E}
\Vertex[x=4,y=0,IdAsLabel]{F}
\Edge[color=blue,label=1,position=left](A)(B)
\Edge[color=blue,label=2,position=above left](B)(C)
\Edge[color=red,label=3,position=above](C)(D)
\Edge[color=blue,label=4,position=above right](D)(E)
\Edge[color=red,label=5,position=right](E)(F)
\Edge[color=red,label=$6_2$,position=below right](A)(D)
\Edge[color=red,label=$7_2$,position=below](A)(F)
\Edge[color=blue,label=8,position=above right](C)(F)
\Edge[color=red,label=9,position=left,bend=15](F)(D)
\Edge[color=blue,label=10,position=above left,bend=5](A)(E)
\Edge[color=red,label=10,position=below right,bend=-5](A)(E)
\end{tikzpicture}
\\
B7(a)
\begin{tikzpicture}[thick,scale=0.8, every node/.style={scale=0.8}]
\Vertex[x=0,y=0,IdAsLabel]{A}
\Vertex[x=0,y=2,IdAsLabel]{B}
\Vertex[x=1,y=4,IdAsLabel]{C}
\Vertex[x=3,y=4,IdAsLabel]{D}
\Vertex[x=4,y=2,IdAsLabel]{E}
\Vertex[x=4,y=0,IdAsLabel]{F}
\Edge[color=blue,label=1,position=left](A)(B)
\Edge[color=red,label=2,position=above left](B)(C)
\Edge[color=blue,label=3,position=above](C)(D)
\Edge[color=blue,label=4,position=above right](D)(E)
\Edge[color=red,label=5,position=right](E)(F)
\Edge[color=blue,label=$6_1$,position=below](A)(F)
\Edge[color=red,label=7,position=above right](C)(F)
\Edge[color=blue,label=8,position=above,bend=5](B)(E)
\Edge[color=red,label=8,position=below,bend=-5](B)(E)
\end{tikzpicture}
B7(b)
\begin{tikzpicture}[thick,scale=0.8, every node/.style={scale=0.8}]
\Vertex[x=0,y=0,IdAsLabel]{A}
\Vertex[x=0,y=2,IdAsLabel]{B}
\Vertex[x=1,y=4,IdAsLabel]{C}
\Vertex[x=3,y=4,IdAsLabel]{D}
\Vertex[x=4,y=2,IdAsLabel]{E}
\Vertex[x=4,y=0,IdAsLabel]{F}
\Vertex[x=2,y=1,IdAsLabel]{G}
\Edge[color=blue,label=1,position=left](A)(B)
\Edge[color=red,label=2,position=above left](B)(C)
\Edge[color=blue,label=3,position=above](C)(D)
\Edge[color=blue,label=4,position=above
right](D)(E)
\Edge[color=red,label=5,position=right](E)(F)
\Edge[color=red,label=$6_2$,position=below](A)(F)
\Edge[color=blue,label=$7_1$,position=above left](A)(G)
\Edge[color=red,label=8,position=above](B)(E)
\Edge[color=red,label=9,position=above right](C)(G)
\Edge[color=blue,label=10,position=above left,bend=5](G)(E)
\Edge[color=red,label=10,position=below right,bend=-5](G)(E)
\end{tikzpicture}
B7(c)
\begin{tikzpicture}[thick,scale=0.8, every node/.style={scale=0.8}]
\Vertex[x=0,y=0,IdAsLabel]{A}
\Vertex[x=0,y=2,IdAsLabel]{B}
\Vertex[x=1,y=4,IdAsLabel]{C}
\Vertex[x=3,y=4,IdAsLabel]{D}
\Vertex[x=4,y=2,IdAsLabel]{E}
\Vertex[x=4,y=0,IdAsLabel]{F}
\Vertex[x=2,y=1,IdAsLabel]{G}
\Edge[color=blue,label=1,position=left](A)(B)
\Edge[color=red,label=2,position=above left](B)(C)
\Edge[color=blue,label=3,position=above](C)(D)
\Edge[color=blue,label=4,position=above right](D)(E)
\Edge[color=red,label=5,position=right](E)(F)
\Edge[color=red,label=$6_2$,position=below](A)(F)
\Edge[color=red,label=$7_2$,position=above left](A)(G)
\Edge[color=blue,label=8,position=below](G)(E)
\Edge[color=red,label=9,position=above left,bend=-15](A)(C)
\Edge[color=blue,label=10,position=below left,bend=5](G)(B)
\Edge[color=red,label=10,position=above right,bend=-5](G)(B)
\end{tikzpicture}
\\
B8(a)
\begin{tikzpicture}[thick,scale=0.8, every node/.style={scale=0.8}]
\Vertex[x=0,y=0,IdAsLabel]{A}
\Vertex[x=0,y=2,IdAsLabel]{B}
\Vertex[x=1,y=4,IdAsLabel]{C}
\Vertex[x=3,y=4,IdAsLabel]{D}
\Vertex[x=4,y=2,IdAsLabel]{E}
\Vertex[x=4,y=0,IdAsLabel]{F}
\Vertex[x=2,y=0,IdAsLabel]{G}
\Edge[color=blue,label=1,position=left](A)(B)
\Edge[color=red,label=2,position=above left](B)(C)
\Edge[color=red,label=3,position=above](C)(D)
\Edge[color=blue,label=4,position=above right](D)(E)
\Edge[color=blue,label=5,position=right](E)(F)
\Edge[color=blue,label=$6_1$,position=above left](A)(D)
\Edge[color=red,label=7,position=below](F)(G)
\Edge[color=red,label=8,position=below left](G)(B)
\Edge[color=blue,label=9,position=above right,bend=5](C)(F)
\Edge[color=red,label=9,position=below left,bend=-5](C)(F)
\end{tikzpicture}
B8(b)
\begin{tikzpicture}[thick,scale=0.8, every node/.style={scale=0.8}]
\Vertex[x=0,y=0,IdAsLabel]{A}
\Vertex[x=0,y=2,IdAsLabel]{B}
\Vertex[x=1,y=4,IdAsLabel]{C}
\Vertex[x=3,y=4,IdAsLabel]{D}
\Vertex[x=4,y=2,IdAsLabel]{E}
\Vertex[x=4,y=0,IdAsLabel]{F}
\Vertex[x=2,y=1,IdAsLabel]{G}
\Edge[color=blue,label=1,position=left](A)(B)
\Edge[color=red,label=2,position=above left](B)(C)
\Edge[color=red,label=3,position=above](C)(D)
\Edge[color=blue,label=4,position=above right](D)(E)
\Edge[color=blue,label=5,position=right](E)(F)
\Edge[color=red,label=$6_2$,position=above left](A)(D)
\Edge[color=blue,label=$7_1$,position=below](A)(F)
\Edge[color=red,label=8,position=below right](G)(D)
\Edge[color=blue,label=9,position=above right,bend=5](B)(G)
\Edge[color=red,label=9,position=below left,bend=-5](B)(G)
\end{tikzpicture}
B8(c)
\begin{tikzpicture}[thick,scale=0.8, every node/.style={scale=0.8}]
\Vertex[x=0,y=0,IdAsLabel]{A}
\Vertex[x=0,y=2,IdAsLabel]{B}
\Vertex[x=1,y=4,IdAsLabel]{C}
\Vertex[x=3,y=4,IdAsLabel]{D}
\Vertex[x=4,y=2,IdAsLabel]{E}
\Vertex[x=4,y=0,IdAsLabel]{F}
\Edge[color=blue,label=1,position=left](A)(B)
\Edge[color=red,label=2,position=above left](B)(C)
\Edge[color=red,label=3,position=above](C)(D)
\Edge[color=blue,label=4,position=above right](D)(E)
\Edge[color=blue,label=5,position=right](E)(F)
\Edge[color=red,label=$6_2$,position=above left](A)(D)
\Edge[color=red,label=$7_2$,position=below](A)(F)
\Edge[color=blue,label=8,position=below left](F)(C)
\Edge[color=red,label=9,position=left,bend=-15](A)(C)
\Edge[color=blue,label=10,position=above left,bend=5](B)(D)
\Edge[color=red,label=10,position=below right,bend=-5](B)(D)
\end{tikzpicture}
\\
B9(a)
\begin{tikzpicture}[thick,scale=0.8, every node/.style={scale=0.8}]
\Vertex[x=0,y=0,IdAsLabel]{A}
\Vertex[x=0,y=2,IdAsLabel]{B}
\Vertex[x=1,y=4,IdAsLabel]{C}
\Vertex[x=3,y=4,IdAsLabel]{D}
\Vertex[x=4,y=2,IdAsLabel]{E}
\Vertex[x=4,y=0,IdAsLabel]{F}
\Edge[color=blue,label=1,position=left](A)(B)
\Edge[color=red,label=2,position=above left](B)(C)
\Edge[color=blue,label=3,position=above](C)(D)
\Edge[color=red,label=4,position=above right](D)(E)
\Edge[color=blue,label=5,position=right](E)(F)
\Edge[color=blue,label=$6_1$,position=above left](B)(D)
\Edge[color=red,label=7,position=above right](C)(E)
\Edge[color=red,label=8,position=below left](C)(F)
\Edge[color=red,label=9,position=above left](A)(E)
\Edge[color=blue,label=10,position=above,bend=5](A)(F)
\Edge[color=red,label=10,position=below,bend=-5](A)(F)
\end{tikzpicture}
B9(b)
\begin{tikzpicture}[thick,scale=0.8, every node/.style={scale=0.8}]
\Vertex[x=0,y=0,IdAsLabel]{A}
\Vertex[x=0,y=2,IdAsLabel]{B}
\Vertex[x=1,y=4,IdAsLabel]{C}
\Vertex[x=3,y=4,IdAsLabel]{D}
\Vertex[x=4,y=2,IdAsLabel]{E}
\Vertex[x=4,y=0,IdAsLabel]{F}
\Edge[color=blue,label=1,position=left](A)(B)
\Edge[color=red,label=2,position=above left](B)(C)
\Edge[color=blue,label=3,position=above](C)(D)
\Edge[color=red,label=4,position=above right](D)(E)
\Edge[color=blue,label=5,position=right](E)(F)
\Edge[color=red,label=$6_2$,position=above left](B)(D)
\Edge[color=blue,label=7,position=above right](C)(E)
\Edge[color=red,label=8,position=below right](A)(D)
\Edge[color=red,label=9,position=above right](B)(F)
\Edge[color=blue,label=10,position=above,bend=5](A)(F)
\Edge[color=red,label=10,position=below,bend=-5](A)(F)
\end{tikzpicture}
\\
R0(a)
\begin{tikzpicture}[thick,scale=0.8, every node/.style={scale=0.8}]
\Vertex[x=0,y=0,IdAsLabel]{A}
\Vertex[x=0,y=2,IdAsLabel]{B}
\Vertex[x=1,y=4,IdAsLabel]{C}
\Vertex[x=3,y=4,IdAsLabel]{D}
\Vertex[x=4,y=2,IdAsLabel]{E}
\Vertex[x=4,y=0,IdAsLabel]{F}
\Edge[color=red,label=1,position=left](A)(B)
\Edge[color=red,label=2,position=above left](B)(C)
\Edge[color=red,label=3,position=above](C)(D)
\Edge[color=red,label=4,position=above right](D)(E)
\Edge[color=red,label=5,position=right](E)(F)
\Edge[color=blue,label=6,position=above left](A)(D)
\Edge[color=blue,label=7,position=below](B)(E)
\Edge[color=blue,label=8,position=above right](F)(C)
\Edge[color=blue,label=$9_1$,position=above right](B)(F)
\Edge[color=red,label=10,position=below right,bend=-15](A)(C)
\Edge[color=blue,label=11,position=above,bend=5](A)(E)
\Edge[color=red,label=11,position=below,bend=-5](A)(E)
\end{tikzpicture}
R0(b)
\begin{tikzpicture}[thick,scale=0.8, every node/.style={scale=0.8}]
\Vertex[x=0,y=0,IdAsLabel]{A}
\Vertex[x=0,y=2,IdAsLabel]{B}
\Vertex[x=1,y=4,IdAsLabel]{C}
\Vertex[x=3,y=4,IdAsLabel]{D}
\Vertex[x=4,y=2,IdAsLabel]{E}
\Vertex[x=4,y=0,IdAsLabel]{F}
\Edge[color=red,label=1,position=left](A)(B)
\Edge[color=red,label=2,position=above left](B)(C)
\Edge[color=red,label=3,position=above](C)(D)
\Edge[color=red,label=4,position=above right](D)(E)
\Edge[color=red,label=5,position=right](E)(F)
\Edge[color=blue,label=6,position=above left](A)(D)
\Edge[color=blue,label=7,position=below](B)(E)
\Edge[color=blue,label=8,position=above right](F)(C)
\Edge[color=red,label=$9_2$,position=above](B)(F)
\Edge[color=blue,label=10,position=below](A)(E)
\Edge[color=blue,label=11,position=above right,bend=5](D)(F)
\Edge[color=red,label=11,position=below left,bend=-5](D)(F)
\end{tikzpicture}
\\
R1(a)
\begin{tikzpicture}[thick,scale=0.8, every node/.style={scale=0.8}]
\Vertex[x=0,y=0,IdAsLabel]{A}
\Vertex[x=0,y=2,IdAsLabel]{B}
\Vertex[x=1,y=4,IdAsLabel]{C}
\Vertex[x=3,y=4,IdAsLabel]{D}
\Vertex[x=4,y=2,IdAsLabel]{E}
\Vertex[x=4,y=0,IdAsLabel]{F}
\Vertex[x=1,y=3,IdAsLabel]{G}
\Edge[color=red,label=1,position=left](A)(B)
\Edge[color=red,label=2,position=above left](B)(C)
\Edge[color=red,label=3,position=above](C)(D)
\Edge[color=red,label=4,position=above right](D)(E)
\Edge[color=blue,label=5,position=right](E)(F)
\Edge[color=blue,label=$6_1$,position=below](A)(F)
\Edge[color=blue,label=7,position=below](B)(E)
\Edge[color=blue,label=8,position=above left](A)(D)
\Edge[color=red,label=9,position=below right](B)(G)
\Edge[color=blue,label=10,position=below right,bend=5](D)(G)
\Edge[color=red,label=10,position=above left,bend=-5](D)(G)
\end{tikzpicture}
R1(b)
\begin{tikzpicture}[thick,scale=0.8, every node/.style={scale=0.8}]
\Vertex[x=0,y=0,IdAsLabel]{A}
\Vertex[x=0,y=2,IdAsLabel]{B}
\Vertex[x=1,y=4,IdAsLabel]{C}
\Vertex[x=3,y=4,IdAsLabel]{D}
\Vertex[x=4,y=2,IdAsLabel]{E}
\Vertex[x=4,y=0,IdAsLabel]{F}
\Edge[color=red,label=1,position=left](A)(B)
\Edge[color=red,label=2,position=above left](B)(C)
\Edge[color=red,label=3,position=above](C)(D)
\Edge[color=red,label=4,position=above right](D)(E)
\Edge[color=blue,label=5,position=right](E)(F)
\Edge[color=red,label=$6_2$,position=below](A)(F)
\Edge[color=blue,label=7,position=below](B)(E)
\Edge[color=blue,label=8,position=above left](A)(D)
\Edge[color=blue,label=9,position=above right](C)(F)
\Edge[color=red,label=10,position=below right](A)(C)
\Edge[color=blue,label=11,position=above left,bend=5](B)(D)
\Edge[color=red,label=11,position=below right,bend=-5](B)(D)
\end{tikzpicture}
\\
R2(a)
\begin{tikzpicture}[thick,scale=0.8, every node/.style={scale=0.8}]
\Vertex[x=0,y=0,IdAsLabel]{A}
\Vertex[x=0,y=2,IdAsLabel]{B}
\Vertex[x=1,y=4,IdAsLabel]{C}
\Vertex[x=3,y=4,IdAsLabel]{D}
\Vertex[x=4,y=2,IdAsLabel]{E}
\Vertex[x=4,y=0,IdAsLabel]{F}
\Vertex[x=3,y=3,IdAsLabel]{G}
\Edge[color=red,label=1,position=left](A)(B)
\Edge[color=red,label=2,position=above left](B)(C)
\Edge[color=red,label=3,position=above](C)(D)
\Edge[color=blue,label=4,position=above right](D)(E)
\Edge[color=red,label=5,position=right](E)(F)
\Edge[color=blue,label=$6_1$,position=below](A)(F)
\Edge[color=blue,label=$7_1$,position=above right](C)(F)
\Edge[color=blue,label=8,position=above left](A)(D)
\Edge[color=red,label=9,position=below](B)(E)
\Edge[color=red,label=10,position=below left](C)(G)
\Edge[color=blue,label=11,position=above right,bend=5](G)(E)
\Edge[color=red,label=11,position=below left,bend=-5](G)(E)
\end{tikzpicture}
R2(b)
\begin{tikzpicture}[thick,scale=0.8, every node/.style={scale=0.8}]
\Vertex[x=0,y=0,IdAsLabel]{A}
\Vertex[x=0,y=2,IdAsLabel]{B}
\Vertex[x=1,y=4,IdAsLabel]{C}
\Vertex[x=3,y=4,IdAsLabel]{D}
\Vertex[x=4,y=2,IdAsLabel]{E}
\Vertex[x=4,y=0,IdAsLabel]{F}
\Vertex[x=2,y=1,IdAsLabel]{G}
\Edge[color=red,label=1,position=left](A)(B)
\Edge[color=red,label=2,position=above left](B)(C)
\Edge[color=red,label=3,position=above](C)(D)
\Edge[color=blue,label=4,position=above right](D)(E)
\Edge[color=red,label=5,position=right](E)(F)
\Edge[color=blue,label=$6_1$,position=below](A)(F)
\Edge[color=red,label=$7_2$,position=above
right](C)(F)
\Edge[color=blue,label=8,position=above left](A)(D)
\Edge[color=blue,label=9,position=below](B)(E)
\Edge[color=red,label=10,position=below left](F)(G)
\Edge[color=blue,label=11,position=above,bend=5](B)(G)
\Edge[color=red,label=11,position=below,bend=-5](B)(G)
\end{tikzpicture}
R2(c)
\begin{tikzpicture}[thick,scale=0.8, every node/.style={scale=0.8}]
\Vertex[x=0,y=0,IdAsLabel]{A}
\Vertex[x=0,y=2,IdAsLabel]{B}
\Vertex[x=1,y=4,IdAsLabel]{C}
\Vertex[x=3,y=4,IdAsLabel]{D}
\Vertex[x=4,y=2,IdAsLabel]{E}
\Vertex[x=4,y=0,IdAsLabel]{F}
\Edge[color=red,label=1,position=left](A)(B)
\Edge[color=red,label=2,position=above left](B)(C)
\Edge[color=red,label=3,position=above](C)(D)
\Edge[color=blue,label=4,position=above right](D)(E)
\Edge[color=red,label=5,position=right](E)(F)
\Edge[color=red,label=$6_2$,position=below](A)(F)
\Edge[color=blue,label=7,position=above right](C)(F)
\Edge[color=blue,label=8,position=above left](A)(D)
\Edge[color=blue,label=9,position=below](B)(E)
\Edge[color=red,label=10,position=below](B)(F)
\Edge[color=blue,label=11,position=above left,bend=-5](A)(C)
\Edge[color=red,label=11,position=right,bend=-15](A)(C)
\end{tikzpicture}
\\
R3(a)
\begin{tikzpicture}[thick,scale=0.8, every node/.style={scale=0.8}]
\Vertex[x=0,y=0,IdAsLabel]{A}
\Vertex[x=0,y=2,IdAsLabel]{B}
\Vertex[x=1,y=4,IdAsLabel]{C}
\Vertex[x=3,y=4,IdAsLabel]{D}
\Vertex[x=4,y=2,IdAsLabel]{E}
\Vertex[x=4,y=0,IdAsLabel]{F}
\Edge[color=red,label=1,position=left](A)(B)
\Edge[color=red,label=2,position=above left](B)(C)
\Edge[color=blue,label=3,position=above](C)(D)
\Edge[color=red,label=4,position=above right](D)(E)
\Edge[color=red,label=5,position=right](E)(F)
\Edge[color=blue,label=$6_1$,position=below](A)(F)
\Edge[color=blue,label=$7_1$,position=above](A)(E)
\Edge[color=blue,label=$8_1$,position=below](B)(F)
\Edge[color=red,label=9,position=below left](C)(E)
\Edge[color=blue,label=10,position=above left,bend=5](B)(D)
\Edge[color=red,label=10,position=below right,bend=-5](B)(D)
\end{tikzpicture}
R3(b)
\begin{tikzpicture}[thick,scale=0.8, every node/.style={scale=0.8}]
\Vertex[x=0,y=0,IdAsLabel]{A}
\Vertex[x=0,y=2,IdAsLabel]{B}
\Vertex[x=1,y=4,IdAsLabel]{C}
\Vertex[x=3,y=4,IdAsLabel]{D}
\Vertex[x=4,y=2,IdAsLabel]{E}
\Vertex[x=4,y=0,IdAsLabel]{F}
\Edge[color=red,label=1,position=left](A)(B)
\Edge[color=red,label=2,position=above left](B)(C)
\Edge[color=blue,label=3,position=above](C)(D)
\Edge[color=red,label=4,position=above right](D)(E)
\Edge[color=red,label=5,position=right](E)(F)
\Edge[color=blue,label=$6_1$,position=below](A)(F)
\Edge[color=blue,label=$7_1$,position=above](A)(E)
\Edge[color=red,label=$8_2$,position=below](B)(F)
\Edge[color=blue,label=9,position=below right](B)(D)
\Edge[color=blue,label=10,position=above right,bend=5](C)(E)
\Edge[color=red,label=10,position=below left,bend=-5](C)(E)
\end{tikzpicture}
R3(c)
\begin{tikzpicture}[thick,scale=0.8, every node/.style={scale=0.8}]
\Vertex[x=0,y=0,IdAsLabel]{A}
\Vertex[x=0,y=2,IdAsLabel]{B}
\Vertex[x=1,y=4,IdAsLabel]{C}
\Vertex[x=3,y=4,IdAsLabel]{D}
\Vertex[x=4,y=2,IdAsLabel]{E}
\Vertex[x=4,y=0,IdAsLabel]{F}
\Edge[color=red,label=1,position=left](A)(B)
\Edge[color=red,label=2,position=above left](B)(C)
\Edge[color=blue,label=3,position=above](C)(D)
\Edge[color=red,label=4,position=above right](D)(E)
\Edge[color=red,label=5,position=right](E)(F)
\Edge[color=blue,label=$6_1$,position=below](A)(F)
\Edge[color=red,label=$7_2$,position=above](A)(E)
\Edge[color=blue,label=8,position=below](B)(F)
\Edge[color=blue,label=9,position=below left](C)(E)
\Edge[color=blue,label=10,position=above left,bend=5](B)(D)
\Edge[color=red,label=10,position=below right,bend=-5](B)(D)
\end{tikzpicture}
\\
R3(d)
\begin{tikzpicture}[thick,scale=0.8, every node/.style={scale=0.8}]
\Vertex[x=0,y=0,IdAsLabel]{A}
\Vertex[x=0,y=2,IdAsLabel]{B}
\Vertex[x=1,y=4,IdAsLabel]{C}
\Vertex[x=3,y=4,IdAsLabel]{D}
\Vertex[x=4,y=2,IdAsLabel]{E}
\Vertex[x=4,y=0,IdAsLabel]{F}
\Edge[color=red,label=1,position=left](A)(B)
\Edge[color=red,label=2,position=above left](B)(C)
\Edge[color=blue,label=3,position=above](C)(D)
\Edge[color=red,label=4,position=above right](D)(E)
\Edge[color=red,label=5,position=right](E)(F)
\Edge[color=red,label=$6_2$,position=below](A)(F)
\Edge[color=blue,label=7,position=below](B)(E)
\Edge[color=blue,label=8,position=above right](C)(F)
\Edge[color=blue,label=9,position=above left](A)(D)
\Edge[color=red,label=10,position=below right](A)(E)
\Edge[color=blue,label=11,position=above right,bend=5](B)(F)
\Edge[color=red,label=11,position=below left,bend=-5](B)(F)
\end{tikzpicture}
\\
R4(a)
\begin{tikzpicture}[thick,scale=0.8, every node/.style={scale=0.8}]
\Vertex[x=0,y=0,IdAsLabel]{A}
\Vertex[x=0,y=2,IdAsLabel]{B}
\Vertex[x=1,y=4,IdAsLabel]{C}
\Vertex[x=3,y=4,IdAsLabel]{D}
\Vertex[x=4,y=2,IdAsLabel]{E}
\Vertex[x=4,y=0,IdAsLabel]{F}
\Vertex[x=2,y=0,IdAsLabel]{G}
\Edge[color=red,label=1,position=left](A)(B)
\Edge[color=red,label=2,position=above left](B)(C)
\Edge[color=red,label=3,position=above](C)(D)
\Edge[color=blue,label=4,position=above right](D)(E)
\Edge[color=blue,label=5,position=right](E)(F)
\Edge[color=blue,label=$6_1$,position=below](A)(G)
\Edge[color=blue,label=7,position=left](A)(D)
\Edge[color=red,label=8,position=right](C)(G)
\Edge[color=red,label=9,position=above right](C)(F)
\Edge[color=red,label=10,position=above right](B)(F)
\Edge[color=blue,label=11,position=above right,bend=5](B)(G)
\Edge[color=red,label=11,position=below left,bend=-5](B)(G)
\end{tikzpicture}
R4(b)
\begin{tikzpicture}[thick,scale=0.8, every node/.style={scale=0.8}]
\Vertex[x=0,y=0,IdAsLabel]{A}
\Vertex[x=0,y=2,IdAsLabel]{B}
\Vertex[x=1,y=4,IdAsLabel]{C}
\Vertex[x=3,y=4,IdAsLabel]{D}
\Vertex[x=4,y=2,IdAsLabel]{E}
\Vertex[x=4,y=0,IdAsLabel]{F}
\Vertex[x=2,y=0,IdAsLabel]{G}
\Edge[color=red,label=1,position=left](A)(B)
\Edge[color=red,label=2,position=above left](B)(C)
\Edge[color=red,label=3,position=above](C)(D)
\Edge[color=blue,label=4,position=above right](D)(E)
\Edge[color=blue,label=5,position=right](E)(F)
\Edge[color=red,label=$6_2$,position=below](A)(G)
\Edge[color=blue,label=7,position=above left](A)(D)
\Edge[color=blue,label=8,position=below left](C)(G)
\Edge[color=red,label=9,position=right](A)(C)
\Edge[color=red,label=10,position=below](G)(F)
\Edge[color=blue,label=11,position=above right,bend=5](C)(F)
\Edge[color=red,label=11,position=below left,bend=-5](C)(F)
\end{tikzpicture}
\\
R5(a)
\begin{tikzpicture}[thick,scale=0.8, every node/.style={scale=0.8}]
\Vertex[x=0,y=0,IdAsLabel]{A}
\Vertex[x=0,y=2,IdAsLabel]{B}
\Vertex[x=1,y=4,IdAsLabel]{C}
\Vertex[x=3,y=4,IdAsLabel]{D}
\Vertex[x=4,y=2,IdAsLabel]{E}
\Vertex[x=4,y=0,IdAsLabel]{F}
\Vertex[x=2,y=1,IdAsLabel]{G}
\Edge[color=blue,label=1,position=left](A)(B)
\Edge[color=red,label=2,position=above left](B)(C)
\Edge[color=red,label=3,position=above](C)(D)
\Edge[color=red,label=4,position=above right](D)(E)
\Edge[color=blue,label=5,position=right](E)(F)
\Edge[color=blue,label=$6_1$,position=below](A)(F)
\Edge[color=blue,label=$7_1$,position=below right](B)(D)
\Edge[color=red,label=8,position=below left](C)(E)
\Edge[color=red,label=9,position=below right](D)(G)
\Edge[color=blue,label=10,position=below right,bend=5](E)(G)
\Edge[color=red,label=10,position=above left,bend=-5](E)(G)
\end{tikzpicture}
R5(b)
\begin{tikzpicture}[thick,scale=0.8, every node/.style={scale=0.8}]
\Vertex[x=0,y=0,IdAsLabel]{A}
\Vertex[x=0,y=2,IdAsLabel]{B}
\Vertex[x=1,y=4,IdAsLabel]{C}
\Vertex[x=3,y=4,IdAsLabel]{D}
\Vertex[x=4,y=2,IdAsLabel]{E}
\Vertex[x=4,y=0,IdAsLabel]{F}
\Vertex[x=2,y=1,IdAsLabel]{G}
\Edge[color=blue,label=1,position=left](A)(B)
\Edge[color=red,label=2,position=above left](B)(C)
\Edge[color=red,label=3,position=above](C)(D)
\Edge[color=red,label=4,position=above right](D)(E)
\Edge[color=blue,label=5,position=right](E)(F)
\Edge[color=blue,label=$6_1$,position=below](A)(F)
\Edge[color=red,label=$7_2$,position=below right](B)(D)
\Edge[color=blue,label=8,position=below left](C)(E)
\Edge[color=red,label=9,position=below left](C)(G)
\Edge[color=blue,label=10,position=below left,bend=-5](B)(G)
\Edge[color=red,label=10,position=above right,bend=5](B)(G)
\end{tikzpicture}
\\
R5(c)
\begin{tikzpicture}[thick,scale=0.8, every node/.style={scale=0.8}]
\Vertex[x=0,y=0,IdAsLabel]{A}
\Vertex[x=0,y=2,IdAsLabel]{B}
\Vertex[x=1,y=4,IdAsLabel]{C}
\Vertex[x=3,y=4,IdAsLabel]{D}
\Vertex[x=4,y=2,IdAsLabel]{E}
\Vertex[x=4,y=0,IdAsLabel]{F}
\Vertex[x=2,y=1,IdAsLabel]{G}
\Edge[color=blue,label=1,position=left](A)(B)
\Edge[color=red,label=2,position=above left](B)(C)
\Edge[color=red,label=3,position=above](C)(D)
\Edge[color=red,label=4,position=above right](D)(E)
\Edge[color=blue,label=5,position=right](E)(F)
\Edge[color=red,label=$6_2$,position=below](A)(F)
\Edge[color=blue,label=$7_1$,position=above left](A)(D)
\Edge[color=blue,label=8,position=below](B)(E)
\Edge[color=red,label=9,position=above right](C)(F)
\Edge[color=red,label=10,position=above left](D)(G)
\Edge[color=blue,label=11,position=below left,bend=5](F)(G)
\Edge[color=red,label=11,position=above right,bend=-5](F)(G)
\end{tikzpicture}
R5(d)
\begin{tikzpicture}[thick,scale=0.8, every node/.style={scale=0.8}]
\Vertex[x=0,y=0,IdAsLabel]{A}
\Vertex[x=0,y=2,IdAsLabel]{B}
\Vertex[x=1,y=4,IdAsLabel]{C}
\Vertex[x=3,y=4,IdAsLabel]{D}
\Vertex[x=4,y=2,IdAsLabel]{E}
\Vertex[x=4,y=0,IdAsLabel]{F}
\Vertex[x=2,y=1,IdAsLabel]{G}
\Edge[color=blue,label=1,position=left](A)(B)
\Edge[color=red,label=2,position=above left](B)(C)
\Edge[color=red,label=3,position=above](C)(D)
\Edge[color=red,label=4,position=above right](D)(E)
\Edge[color=blue,label=5,position=right](E)(F)
\Edge[color=red,label=$6_2$,position=below](A)(F)
\Edge[color=red,label=$7_2$,position=above left](A)(D)
\Edge[color=blue,label=8,position=below](B)(E)
\Edge[color=blue,label=9,position=above right](C)(F)
\Edge[color=red,label=10,position=above right](C)(G)
\Edge[color=blue,label=11,position=below right,bend=-5](A)(G)
\Edge[color=red,label=11,position=above left,bend=5](A)(G)
\end{tikzpicture}
\\
R6(a)
\begin{tikzpicture}[thick,scale=0.8, every node/.style={scale=0.8}]
\Vertex[x=0,y=0,IdAsLabel]{A}
\Vertex[x=0,y=2,IdAsLabel]{B}
\Vertex[x=1,y=4,IdAsLabel]{C}
\Vertex[x=3,y=4,IdAsLabel]{D}
\Vertex[x=4,y=2,IdAsLabel]{E}
\Vertex[x=4,y=0,IdAsLabel]{F}
\Edge[color=red,label=1,position=left](A)(B)
\Edge[color=red,label=2,position=above left](B)(C)
\Edge[color=blue,label=3,position=above](C)(D)
\Edge[color=red,label=4,position=above right](D)(E)
\Edge[color=blue,label=5,position=right](E)(F)
\Edge[color=blue,label=$6_1$,position=below](A)(F)
\Edge[color
=blue,label=$7_1$,position=below](B)(E)
\Edge[color=red,label=8,position=right](A)(C)
\Edge[color=red,label=9,position=above left](A)(D)
\Edge[color=blue,label=10,position=above left,bend=5](B)(D)
\Edge[color=red,label=10,position=below right,bend=-5](B)(D)
\end{tikzpicture}
R6(b)
\begin{tikzpicture}[thick,scale=0.8, every node/.style={scale=0.8}]
\Vertex[x=0,y=0,IdAsLabel]{A}
\Vertex[x=0,y=2,IdAsLabel]{B}
\Vertex[x=1,y=4,IdAsLabel]{C}
\Vertex[x=3,y=4,IdAsLabel]{D}
\Vertex[x=4,y=2,IdAsLabel]{E}
\Vertex[x=4,y=0,IdAsLabel]{F}
\Vertex[x=1,y=2.5,IdAsLabel]{G}
\Edge[color=red,label=1,position=left](A)(B)
\Edge[color=red,label=2,position=above left](B)(C)
\Edge[color=blue,label=3,position=above](C)(D)
\Edge[color=red,label=4,position=above right](D)(E)
\Edge[color=blue,label=5,position=right](E)(F)
\Edge[color=blue,label=$6_1$,position=below](A)(F)
\Edge[color=red,label=$7_2$,position=below](B)(E)
\Edge[color=blue,label=8,position=above left](A)(D)
\Edge[color=red,label=9,position=above](E)(G)
\Edge[color=blue,label=10,position=right,bend=5](C)(G)
\Edge[color=red,label=10,position=left,bend=-5](C)(G)
\end{tikzpicture}
\\
R6(c)
\begin{tikzpicture}[thick,scale=0.8, every node/.style={scale=0.8}]
\Vertex[x=0,y=0,IdAsLabel]{A}
\Vertex[x=0,y=2,IdAsLabel]{B}
\Vertex[x=1,y=4,IdAsLabel]{C}
\Vertex[x=3,y=4,IdAsLabel]{D}
\Vertex[x=4,y=2,IdAsLabel]{E}
\Vertex[x=4,y=0,IdAsLabel]{F}
\Vertex[x=1,y=3,IdAsLabel]{G}
\Edge[color=red,label=1,position=left](A)(B)
\Edge[color=red,label=2,position=above left](B)(C)
\Edge[color=blue,label=3,position=above](C)(D)
\Edge[color=red,label=4,position=above right](D)(E)
\Edge[color=blue,label=5,position=right](E)(F)
\Edge[color=red,label=$6_2$,position=below](A)(F)
\Edge[color=blue,label=$7_1$,position=below](B)(E)
\Edge[color=blue,label=8,position=above right](C)(F)
\Edge[color=red,label=9,position=above left](A)(D)
\Edge[color=red,label=10,position=above left](G)(D)
\Edge[color=blue,label=11,position=above left,bend=5](B)(G)
\Edge[color=red,label=11,position=below right,bend=-5](B)(G)
\end{tikzpicture}
R6(d)
\begin{tikzpicture}[thick,scale=0.8, every node/.style={scale=0.8}]
\Vertex[x=0,y=0,IdAsLabel]{A}
\Vertex[x=0,y=2,IdAsLabel]{B}
\Vertex[x=1,y=4,IdAsLabel]{C}
\Vertex[x=3,y=4,IdAsLabel]{D}
\Vertex[x=4,y=2,IdAsLabel]{E}
\Vertex[x=4,y=0,IdAsLabel]{F}
\Vertex[x=2,y=1,IdAsLabel]{G}
\Edge[color=red,label=1,position=left](A)(B)
\Edge[color=red,label=2,position=above left](B)(C)
\Edge[color=blue,label=3,position=above](C)(D)
\Edge[color=red,label=4,position=above right](D)(E)
\Edge[color=blue,label=5,position=right](E)(F)
\Edge[color=red,label=$6_2$,position=below](A)(F)
\Edge[color=red,label=$7_2$,position=below](B)(E)
\Edge[color=blue,label=8,position=above right](C)(F)
\Edge[color=blue,label=9,position=above left](A)(D)
\Edge[color=red,label=10,position=below right](A)(G)
\Edge[color=blue,label=11,position=below right,bend=5](E)(G)
\Edge[color=red,label=11,position=above left,bend=-5](E)(G)
\end{tikzpicture}
\\
R7(a)
\begin{tikzpicture}[thick,scale=0.8, every node/.style={scale=0.8}]
\Vertex[x=0,y=0,IdAsLabel]{A}
\Vertex[x=0,y=2,IdAsLabel]{B}
\Vertex[x=1,y=4,IdAsLabel]{C}
\Vertex[x=3,y=4,IdAsLabel]{D}
\Vertex[x=4,y=2,IdAsLabel]{E}
\Vertex[x=4,y=0,IdAsLabel]{F}
\Edge[color=red,label=1,position=left](A)(B)
\Edge[color=blue,label=2,position=above left](B)(C)
\Edge[color=red,label=3,position=above](C)(D)
\Edge[color=red,label=4,position=above right](D)(E)
\Edge[color=blue,label=5,position=right](E)(F)
\Edge[color=blue,label=$6_1$,position=below right](B)(D)
\Edge[color=blue,label=$7_1$,position=below right](A)(C)
\Edge[color=red,label=8,position=above left](A)(E)
\Edge[color=red,label=9,position=left](D)(F)
\Edge[color=blue,label=10,position=above,bend=5](A)(F)
\Edge[color=red,label=10,position=below,bend=-5](A)(F)
\end{tikzpicture}
R7(b)
\begin{tikzpicture}[thick,scale=0.8, every node/.style={scale=0.8}]
\Vertex[x=0,y=0,IdAsLabel]{A}
\Vertex[x=0,y=2,IdAsLabel]{B}
\Vertex[x=1,y=4,IdAsLabel]{C}
\Vertex[x=3,y=4,IdAsLabel]{D}
\Vertex[x=4,y=2,IdAsLabel]{E}
\Vertex[x=4,y=0,IdAsLabel]{F}
\Edge[color=red,label=1,position=left](A)(B)
\Edge[color=blue,label=2,position=above left](B)(C)
\Edge[color=red,label=3,position=above](C)(D)
\Edge[color=red,label=4,position=above right](D)(E)
\Edge[color=blue,label=5,position=right](E)(F)
\Edge[color=blue,label=$6_1$,position=below right](B)(D)
\Edge[color=red,label=$7_2$,position=below right](A)(C)
\Edge[color=blue,label=8,position=above left](A)(E)
\Edge[color=red,label=9,position=above left](A)(D)
\Edge[color=red,label=10,position=above right](C)(F)
\Edge[color=blue,label=11,position=right,bend=-5](D)(F)
\Edge[color=red,label=11,position=left,bend=-15](D)(F)
\end{tikzpicture}
R7(c)
\begin{tikzpicture}[thick,scale=0.8, every node/.style={scale=0.8}]
\Vertex[x=0,y=0,IdAsLabel]{A}
\Vertex[x=0,y=2,IdAsLabel]{B}
\Vertex[x=1,y=4,IdAsLabel]{C}
\Vertex[x=3,y=4,IdAsLabel]{D}
\Vertex[x=4,y=2,IdAsLabel]{E}
\Vertex[x=4,y=0,IdAsLabel]{F}
\Vertex[x=2,y=0,IdAsLabel]{G}
\Edge[color=red,label=1,position=left](A)(B)
\Edge[color=blue,label=2,position=above left](B)(C)
\Edge[color=red,label=3,position=above](C)(D)
\Edge[color=red,label=4,position=above right](D)(E)
\Edge[color=blue,label=5,position=right](E)(F)
\Edge[color=red,label=$6_2$,position=below right](B)(D)
\Edge[color=blue,label=7,position=below right](A)(C)
\Edge[color=blue,label=8,position=above left](A)(E)
\Edge[color=red,label=9,position=below](G)(F)
\Edge[color=red,label=10,position=below left](B)(G)
\Edge[color=blue,label=11,position=right,bend=-5](D)(F)
\Edge[color=red,label=11,position=left,bend=-15](D)(F)
\end{tikzpicture}
\\
R8(a)
\begin{tikzpicture}[thick,scale=0.8, every node/.style={scale=0.8}]
\Vertex[x=0,y=0,IdAsLabel]{A}
\Vertex[x=0,y=2,IdAsLabel]{B}
\Vertex[x=1,y=4,IdAsLabel]{C}
\Vertex[x=3,y=4,IdAsLabel]{D}
\Vertex[x=4,y=2,IdAsLabel]{E}
\Vertex[x=4,y=0,IdAsLabel]{F}
\Vertex[x=2,y=0,IdAsLabel]{G}
\Edge[color=red,label=1,position=left](A)(B)
\Edge[color=blue,label=2,position=above left](B)(C)
\Edge[color=blue,label=3,position=above](C)(D)
\Edge[color=red,label=4,position=above right](D)(E)
\Edge[color=red,label=5,position=right](E)(F)
\Edge[color=blue,label=$6_1$,position=below](B)(E)
\Edge[color=blue,label=$7_1$,position=below](G)(F)
\Edge[color=red,label=8,position=above left](G)(E)
\Edge[color=red,label=9,position=below left](D)(F)
\Edge[color=blue,label=10,position=above left,bend=5](G)(D)
\Edge[color=red,label=10,position=below right,bend=-5](G)(D)
\end{tikzpicture}
R8(b)
\begin{tikzpicture}[thick,scale=0.8, every node/.style={scale=0.8}]
\Vertex[x=0,y=0,IdAsLabel]{A}
\Vertex[x=0,y=2,IdAsLabel]{B}
\Vertex[x=1,y=4,IdAsLabel]{C}
\Vertex[x=3,y=4,IdAsLabel]{D}
\Vertex[x=4,y=2,IdAsLabel]{E}
\Vertex[x=4,y=0,IdAsLabel]{F}
\Vertex[x=2,y=0,IdAsLabel]{G}
\Edge[color=red,label=1,position=left](A)(B)
\Edge[color=blue,label=2,position=above left](B)(C)
\Edge[color=blue,label=3,position=above](C)(D)
\Edge[color=red,label=4,position=above right](D)(E)
\Edge[color=red,label=5,position=right](E)(F)
\Edge[color=blue,label=$6_1$,position=below](B)(E)
\Edge[color=red,label=$7_2$,position=below](G)(F)
\Edge[color=blue,label=8,position=above left](G)(D)
\Edge[color=red,label=9,position=above left](A)(E)
\Edge[color=blue,label=10,position=above,bend=5](A)(G)
\Edge[color=red,label=10,position=below,bend=-5](A)(G)
\end{tikzpicture}
R8(c)
\begin{tikzpicture}[thick,scale=0.8, every node/.style={scale=0.8}]
\Vertex[x=0,y=0,IdAsLabel]{A}
\Vertex[x=0,y=2,IdAsLabel]{B}
\Vertex[x=1,y=4,IdAsLabel]{C}
\Vertex[x=3,y=4,IdAsLabel]{D}
\Vertex[x=4,y=2,IdAsLabel]{E}
\Vertex[x=4,y=0,IdAsLabel]{F}
\Vertex[x=2,y=1,IdAsLabel]{G}
\Edge[color=red,label=1,position=left](A)(B)
\Edge[color=blue,label=2,position=above left](B)(C)
\Edge[color=blue,label=3,position=above](C)(D)
\Edge[color=red,label=4,position=above right](D)(E)
\Edge[color=red,label=5,position=right](E)(F)
\Edge[color=red,label=$6_2$,position=below](B)(E)
\Edge[color=blue,label=7,position=above left](A)(D)
\Edge[color=blue,label=8,position=below](A)(F)
\Edge[color=red,label=9,position=below right](B)(G)
\Edge[color=blue,label=10,position=below left,bend=5](F)(G)
\Edge[color=red,label=10,position=above right,bend=-5](F)(G)
\end{tikzpicture}
\\
R9(a)
\begin{tikzpicture}[thick,scale=0.8, every node/.style={scale=0.8}]
\Vertex[x=0,y=0,IdAsLabel]{A}
\Vertex[x=0,y=2,IdAsLabel]{B}
\Vertex[x=1,y=4,IdAsLabel]{C}
\Vertex[x=3,y=4,IdAsLabel]{D}
\Vertex[x=4,y=2,IdAsLabel]{E}
\Vertex[x=4,y=0,IdAsLabel]{F}
\Edge[color=red,label=1,position=left](A)(B)
\Edge[color=blue,label=2,position=above left](B)(C)
\Edge[color=red,label=3,position=above](C)(D)
\Edge[color=blue,label=4,position=above right](D)(E)
\Edge[color=red,label=5,position=right](E)(F)
\Edge[color=blue,label=$6_1$,position=below](A)(F)
\Edge[color=blue,label=$7_1$,position=above](A)(E)
\Edge[color=red,label=8,position=below](B)(F)
\Edge[color=red,label=9,position=above right](C)(F)
\Edge[color=blue,label=10,position=below right,bend=-5](B)(D)
\Edge[color=red,label=10,position=above left,bend=5](B)(D)
\end{tikzpicture}
R9(b)
\begin{tikzpicture}[thick,scale=0.8, every node/.style={scale=0.8}]
\Vertex[x=0,y=0,IdAsLabel]{A}
\Vertex[x=0,y=2,IdAsLabel]{B}
\Vertex[x=1,y=4,IdAsLabel]{C}
\Vertex[x=3,y=4,IdAsLabel]{D}
\Vertex[x=4,y=2,IdAsLabel]{E}
\Vertex[x=4,y=0,IdAsLabel]{F}
\Edge[color=red,label=1,position=left](A)(B)
\Edge[color=blue,label=2,position=above left](B)(C)
\Edge[color=red,label=3,position=above](C)(D)
\Edge[color=blue,label=4,position=above right](D)(E)
\Edge[color=red,label=5,position=right](E)(F)
\Edge[color=blue,label=$6_1$,position=below](A)(F)
\Edge[color=red,label=$7_2$,position=above](A)(E)
\Edge[color=blue,label=8,position=below](B)(F)
\Edge[color=red,label=9,position=above left](A)(D)
\Edge[color=blue,label=10,position=below left,bend=-5](C)(E)
\Edge[color=red,label=10,position=above right,bend=5](C)(E)
\end{tikzpicture}
\\
R9(c)
\begin{tikzpicture}[thick,scale=0.8, every node/.style={scale=0.8}]
\Vertex[x=0,y=0,IdAsLabel]{A}
\Vertex[x=0,y=2,IdAsLabel]{B}
\Vertex[x=1,y=4,IdAsLabel]{C}
\Vertex[x=3,y=4,IdAsLabel]{D}
\Vertex[x=4,y=2,IdAsLabel]{E}
\Vertex[x=4,y=0,IdAsLabel]{F}
\Vertex[x=2,y=1,IdAsLabel]{G}
\Edge[color=red,label=1,position=left](A)(B)
\Edge[color=blue,label=2,position=above left](B)(C)
\Edge[color=red,label=3,position=above](C)(D)
\Edge[color=blue,label=4,position=above right](D)(E)
\Edge[color=red,label=5,position=right](E)(F)
\Edge[color=red,label=$6_2$,position=below](A)(F)
\Edge[color=blue,label=$7_1$,position=above
left](A)(D)
\Edge[color=blue,label=8,position=below](B)(E)
\Edge[color=red,label=9,position=above right](C)(F)
\Edge[color=red,label=10,position=above right](C)(G)
\Edge[color=blue,label=11,position=below right,bend=-5](A)(G)
\Edge[color=red,label=11,position=above left,bend=5](A)(G)
\end{tikzpicture}
R9(d)
\begin{tikzpicture}[thick,scale=0.8, every node/.style={scale=0.8}]
\Vertex[x=0,y=0,IdAsLabel]{A}
\Vertex[x=0,y=2,IdAsLabel]{B}
\Vertex[x=1,y=4,IdAsLabel]{C}
\Vertex[x=3,y=4,IdAsLabel]{D}
\Vertex[x=4,y=2,IdAsLabel]{E}
\Vertex[x=4,y=0,IdAsLabel]{F}
\Vertex[x=2,y=1,IdAsLabel]{G}
\Edge[color=red,label=1,position=left](A)(B)
\Edge[color=blue,label=2,position=above left](B)(C)
\Edge[color=red,label=3,position=above](C)(D)
\Edge[color=blue,label=4,position=above right](D)(E)
\Edge[color=red,label=5,position=right](E)(F)
\Edge[color=red,label=$6_2$,position=below](A)(F)
\Edge[color=red,label=$7_2$,position=above left](A)(D)
\Edge[color=blue,label=8,position=below](B)(E)
\Edge[color=blue,label=9,position=above right](C)(F)
\Edge[color=red,label=10,position=above left](D)(G)
\Edge[color=blue,label=11,position=below left,bend=5](F)(G)
\Edge[color=red,label=11,position=above right,bend=-5](F)(G)
\end{tikzpicture}
\newline\newline \indent
The edges marked with the numbers $1$ to $5$ correspond to five consecutive edges of the path $P_6$ from the appropriate row in the table in Figure \ref{1}. Their colors were imposed by considering all possible cases of coloring these first five edges. At this stage, we have six vertices of the path $P_6$ marked with successive (uppercase) letters of the alphabet $(A,B,C,D,E,F)$.\newline \indent
The edges marked with numbers from $6$ to $11$ correspond to the next edges selected by Builder. Their colors have been forced by the fact that Painter wants to delay Builder for as many rounds as possible. So red color of the edge means that if Painter had colored this edge blue, he would have created blue $P_6$. Similarly, blue color of the edge means that if Painter had colored this edge red, he would have created red $C_4$. If a given subcase required it, we added a seventh vertex (marked by $G$).\newline \indent
A number with a subscript ($N_1$ and $N_2$) means that Painter had a choice on the color of that edge in this round, which led to branching into subcases. We will consider a maximum of four subcases, which we denote by successive (lowercase) letters of the alphabet $(a,b,c,d)$. Whereas (for the sake of order), we always assume that the $N_1$ edge is blue and the $N_2$ edge is red. We have a total of $50$ subcases to consider.\newline \indent
Builder forces a blue path $P_6$ or a red cycle $C_4$ over the next six rounds, as shown in the graphs above. The last edge is drawn in two colors. Note that after Builder selects this edge, Painter, by coloring it blue, will create a blue $P_6$, and by coloring it red, he will create a red $C_4$. This prevents Painter from making another non-losing move.\newline\newline \indent
It would be too long and impractical to discuss each subcase in detail. Therefore, we will now make some general remarks and comments about the above proof.\newline \indent
First, let us note that some subcases are symmetrical. This is frequent because, given two consecutively selected edges, Painter has a choice as to the color of the first of these edges, but this choice forces the second of these edges to be colored with the second color. So, regardless of the color choice of the first edge, after $2$ rounds we have the same situation. Examples of symmetry are: $B5(a-b)$; $B9(a-b)$; $R3(b-c)$; $R5(a-b)$; $R5(c-d)$; $R9(a-b)$; $R9(c-d)$. In these situations, we use the same strategy for Builder in both subcases.\newline \indent
Next, let us note that some subcases can be reduced to others. This is often due to the fact that the sixth edge is the edge between vertices $A$ and $F$, so the first $6$ edges form the cycle $C_6$. So if in a given cycle the arrangement of colors is the same (with the accuracy of a rotation) as in another subcase, then the considered subcase is reduced to it. Examples of transitions to the subcase already considered are: $R3(d)$ to $R2(c)$ to $R1(b)$; $R6(c)$ to $R2(a)$; $R6(d)$ to $R2(b)$; $R8(c)$ to $R6(b)$. In these situations, we use the same strategy for Builder in both subcases.\newline \indent
One can also consider whether we have applied the optimal strategy in each of the cases considered. Note, however, that we don't really need optimal strategies. Since we know by Proposition \ref{bounds} that $11\leq\tilde{r}(C_4,P_6)$, we only need to point out strategies for Builder to win in $11$ rounds. Another thing to consider is that we always start with $P_6$. However, as above, it suffices to show that Builder can win in $11$ rounds.\newline \indent
Finally, we can describe one case in detail. Let's choose the case $R3$, because since it is the only case where there are three edges where Painter had a choice, it may seem the most complicated. In addition, we can discuss the above comments on this example. We will therefore analyze each of the four subcases in detail.
\begin{itemize}
	\item $R3(a)$ The edges marked with the numbers $1$ to $5$ correspond to five consecutive edges of the path $P_6$ from the appropriate row $28$ in the table in Figure \ref{1}. Their colors were imposed by considering all possible cases of coloring these first five edges. In our case, we have a color scheme: red, red, blue, red, red. The sixth edge chosen by Builder is the edge between vertices $A$ and $F$. In short we will write: edge $(A,F)$. In this case, we have a number with a subscript ($6_1$ and $6_2$), which means that Painter had a choice about the color of the edge in this round, which led to branching into subcases. By convention, we assume that $6_1$ is blue and $6_2$ is red. A similar situation occurs for the seventh edge $(A,E)$ and the eighth edge $(B,F)$. In our subcase $R3(a)$ all of these three edges are colored blue. Now note that the ninth edge $(C,E)$ must be red, otherwise Painter would create the blue $P_6:(B,F,A,E,C,D)$. The last, tenth edge $(B,D)$ is drawn in two colors. Note that after Builder selects this edge, Painter, by coloring it blue, will create blue $P_6:(C,D,B,F,A,E)$, and by coloring it red, he will create red $C_4:(B,C,E,D)$. This makes it impossible to make another non-losing move.
	\item $R3(b)$ In this subcase, the sixth and seventh edges are colored blue (as in $R3(a)$), but the eighth edges are colored red. Now notice that the ninth edge $(B,D)$ must be blue, otherwise Painter would create red $C_4:(B,D,E,F)$. The last, tenth edge $(C,E)$ is drawn in two colors. Note that after Builder selects this edge, Painter, by coloring it blue, will create blue $P_6:(B,D,C,E,A,F)$, and by coloring it red, he will create red $C_4:(B,C,E,F)$. This makes it impossible to make another non-losing move.
	\item $R3(c)$ In this subcase, the sixth edge is colored blue (as in $R3(b)$), but the seventh edge is colored red. Note that the eighth edge $(B,F)$ must be blue, otherwise Painter would create red $C_4:(A,B,F,E)$. Now note that this subcase can be solved by symmetry with the subcase $R3(b)$. Considering the two consecutively selected edges (seventh and eighth), we see that after $2$ rounds we actually have the same situation. So we use the same strategy for Builder as in the $R3(b)$ subcase.
	\item $R3(d)$ Note that this subcase differs significantly from the others because the sixth edge $(A,F)$ is colored red, not blue as in the previous three subcases. Note, however, that this subcase can be reduced to $R1(b)$ and $R2(c)$. This is because the sixth edge is the edge between the vertices $A$ and $F$, so the first $6$ edges form the cycle $C_6$. So it suffices to note that in our cycle the arrangement of colors is the same as in the two previous subcases. So the subcase under consideration boils down to them. In these situations, we use the same strategy for Builder as in those subcases. However, since we are not discussing all the subcases one by one, we need to discuss this subcase now. So note that the seventh edge $(B,E)$ must be blue, otherwise Painter would create red $C_4:(A,B,E,F)$. Likewise, the eighth edge $(C,F)$ must be blue, otherwise Painter would create red $C_4:(A,B,C,F)$. Similarly, the ninth edge $(A,D)$ must be blue, otherwise Painter would create red $C_4:(A,D,E,F)$. Finally, note that the tenth edge $(A,E)$ must be red, otherwise Painter would create the blue $P_6:(B,E,A,D,C,F)$. The last, eleventh edge $(B,F)$ is drawn in two colors. Note that after Builder selects this edge, Painter, by coloring it blue, will create blue $P_6:(A,D,C,F,B,E)$, and by coloring it red, he will create red $C_4:(A,B,F,E)$. This makes it impossible to make another non-losing move.
\end{itemize}

\bibliographystyle{amsplain}

\end{document}